\documentclass[12pt,reqno]{amsart}
\usepackage{color}
\usepackage{amssymb}

\newtheorem{theorem}{Theorem}[section]
\newtheorem{lemma}[theorem]{Lemma}
\newtheorem{proposition}[theorem]{Proposition}

\theoremstyle{definition}

\newtheorem{remark}[theorem]{Remark}

\numberwithin{equation}{section}

\begin{document}

\baselineskip=15.5pt

\title{Ideals and conjugacy classes in solvable Lie algebras}

\author[S. Ali]{Sajid Ali}

\address{College of Interdisciplinary Studies, Zayed University, 19282 Dubai, UAE}

\email{sajid.ali@zu.ac.ae}

\author[H. Azad]{Hassan Azad}

\address{Abdus Salam School of Mathematical Sciences, GCU, Lahore 54600, Pakistan}

\email{hassan.azad@sms.edu.pk}

\author[I. Biswas]{Indranil Biswas}

\address{School of Mathematics, Tata Institute of Fundamental
Research, Homi Bhabha Road, Mumbai 400005, India}

\email{indranil@math.tifr.res.in}

\author[F. M. Mahomed]{Fazal M. Mahomed}

\address{DSI-NRF Centre of Excellence in Mathematical and Statistical Sciences, School of Computer Science and Applied
Mathematics, University of the Witwatersrand, Johannesburg, Wits 2050, South Africa}

\email{Fazal.Mahomed@wits.edu.za}

\subjclass[2000]{34A26\, 37C10\, 57R30. }

\keywords{Solvable Lie algebra, conjugacy class, ideal}

\date{}

\begin{abstract}
A constructive procedure is given to determine all ideals of a
finite dimensional solvable Lie algebra. This is used in determining
all conjugacy classes of subalgebras of a given finite dimensional solvable Lie algebra.
\end{abstract}

\maketitle

\tableofcontents

\section{Introduction}\label{se1}

All algebras considered in this paper are finite dimensional.
Conjugacy classes of subalgebras of a Lie algebra are considered in the thesis of Ibragimov \cite{ib1} and this has been used to construct
solutions of differential equations of physical interest in
\cite{ib1,ib2,ib3}. The method is explained briefly in \cite[p.~5--7]{ib2}. It involves difficult algebraic
computations and it is not algorithmic. This suggests the problem of devising an algorithmic procedure to find conjugacy classes of
subalgebras. Such an algorithm was developed by Amata and Oliveri in \cite{AO}.
Some of the examples considered in this paper were also considered in the paper of
Patera and Winternitz \cite{PW}. An earlier paper on the same subject is \cite{kotz}.

In this note, we take a different approach: namely we give an algorithmic procedure to find all ideals of a solvable Lie algebra and we use
this knowledge to determine conjugacy classes of all subalgebras of a given solvable Lie algebra. The crucial case turns out
to be one dimensional subalgebras.
The reason is that if $H$ is a subalgebra of a solvable algebra $g$, then the one dimensional extensions of $H$ live in $N(H)/H$ and any
solvable Lie algebra can be obtained by a sequence of one dimensional extensions.

The connection of conjugacy classes with ideals comes from the following considerations: every 
element generates an ideal, which could be the full algebra. Given an element $X$, we say that
$X$ is of shape ${\mathcal I}$ if ${\mathcal I}$ is the ideal of smallest dimension containing $X$.
Therefore, a list of all ideals determines the form of conjugacy 
classes of one dimensional subalgebras. Factorization of the adjoint group related to the ideal 
simplifies calculations considerably.

For semisimple Lie algebras, the classes of one dimensional subalgebras need a description of all maximal solvable Lie algebras. This can be
done, adapting the arguments in \cite{has1}.

Regarding semisimple algebras and their semisimple subalgebras, different ideas are needed that are discussed in detail by de Graff
\cite{graaf}, and Minchenko \cite{min}. Ideals of solvable Lie algebras can also be obtained using \cite{has2}. The main difference is that
the algorithm of this paper involves eigenvalue computations in relatively low dimensions. However we have indicated briefly in Example 4.4
how a knowledge of conjugacy classes of solvable algebras is helpful in determining all conjugacy classes in nonsolvable algebras.

This paper is organized as follows: In Section 2, we give an algorithm to find all ideals of a 
solvable Lie algebra. In Section 3, we implement the algorithm and use it to find conjugacy 
classes of all subalgebras in certain low dimensional Lie algebras. Section 4 contains some basic 
examples.

\section{Ideals of Solvable Lie Algebras}

The algorithm to find all the ideals of solvable Lie algebra is based on the following proposition:

\begin{proposition}\label{prop1}
Let $g$ be a complex solvable Lie algebra. Then a $1$-dimensional ideal of it is either in the center $Z(g)$ or it is spanned by an element $X$
lying in the center $Z(g')$ of $g'\,:=\,[g,\, g]$ and satisfying the condition that $X$
is an eigenvector for the abelian algebra $g/\langle g',\,Z(g)\rangle $ operating on $Z(g')$.
\end{proposition}

First notice that if an ideal ${\mathcal I}\, \subset\, {\mathfrak g}$ centralizes
an element $X$, then we have an well-defined action of ${\mathfrak g}/{\mathcal I}$ on $X$, namely
$\text{ad}(\overline{Y})(X)\,=\, [Y,\, X]$, where $\overline{Y}\,=\, Y+{\mathcal I}$.

\begin{proof}[{Proof of Proposition \ref{prop1}}]
Let $\langle X\rangle $ be a $1$-dimensional ideal of $g$ which is not contained in the center
$Z(g)$ of $g$. Thus there is an element $Y$ in $g$
with $[Y,\,X]\,=\,\lambda\cdot X$ with $\lambda \,\neq\, 0$. Hence, $X$ is in $g'$. Now if $U \,\in\, g'$, then $[U,\,X]\,=\,0$ because ${\rm ad}(U)$
is nilpotent. Hence, $X \,\in\, Z(g')$. Therefore, $\langle g',\, Z(g)\rangle $ operates trivially on $X$ and
consequently $X$ is a common eigenvector for the abelian algebra $g/\langle g',\, Z(g)\rangle $. The converse is obvious.
\end{proof}

\subsection{Algorithm for ideals of a complex solvable Lie algebra}
In this subsection, we give an algorithm for computing all ideals of a complex solvable Lie algebra
as well as expressing them as a sequence of one--dimensional extensions.

Let $g$ be a complex solvable Lie algebra.
By Lie's theorem (see \cite{HN}) we know that every ideal $I$ of $g$ has a series of ideals
$$0\,=\, I_0\,\subset\, I_{1} \,\subset\, I_{2} \,\subset\, I_{3} \,\subset\, \cdots \,\subset \,I_{d}\,=\,I$$ such that each
successive quotient is of dimension one.

\textbf{Step 1.} Extending a basis of
$\langle g',\, Z(g)\rangle $ to a basis of $g$, say $$g\,=\,\langle g', \,Z(g)\rangle \oplus \langle e_{1}\rangle
\oplus \dots \oplus \langle e_{r}\rangle ,$$ we see that the projections
$\widetilde{e}_{1},\, \cdots,\, \widetilde{e}_{r}$ form a basis of $\widetilde{g} \,= \,g/\langle g',\, Z(g)\rangle $ with the algebra
$\widetilde{g}$ being abelian and operating on $Z(g')$. Thus any common eigenvector of $\widetilde{e}_{1},\, \cdots, \,\widetilde{e}_{r}$
on $Z(g')$ gives a $1$-dimensional ideal of $g$. Any non-zero subspace of $Z(g)$ is evidently an ideal of $g$.

\textbf{Step 2.} Assume inductively that all $d$-dimensional ideals have been constructed. Let $J$ be a $(d+1)$--dimensional ideal. Then $J$
has an ideal $I$ of codimension $1$. Thus, $J/I \,\subset\, g/I$ is a $1$--dimensional ideal of $g/I$. Since $I$ is known, we can repeat Step 2 on
$g/I$ to construct $J/I$ and $J$ itself. This completes the procedure for finding all ideals of a complex solvable Lie algebra.

\subsection{Algorithm for ideals of a real solvable Lie algebra}

In this subsection, we give an algorithm for computing all ideals of a real Lie algebra
and writing them as a sequence of extensions of codimension at most two.

Let $g$ be a real solvable Lie algebra, and let $g^{\mathbb{C}}\,=\,g\oplus\sqrt{-1}g$ be its complexification. If $I$ is an ideal of $g$, then
$I^{\mathbb{C}}\,=\, I\otimes {\mathbb C}$ is an ideal of $g^{\mathbb{C}}$, say of dimension $d$. By the procedure in Section 2.1, the
ideal $I^{\mathbb{C}}$ is in a chain of ideals of $g^{\mathbb{C}}$, say $$0\,=\,I_{0}\,\subset\, I_{1}\,\subset\, \cdots\,\subset\, I_{d}
\,=\,I^{\mathbb{C}}$$ with $\dim_{\mathbb{C}}(I_{k}/I_{k-1})\,=\,1$. Therefore, $$0\,\subset\, I_{1} \cap g\,\subset \,\cdots
\, \subset\, I_{d} \cap g \,=\, I$$ is a sequence of ideals of $g$
with $\dim_{\mathbb{R}}(I_{k}\cap g/I_{k-1})\cap g) \,\leq\, 2$. This gives a constructive procedure for finding all
ideals of a real solvable Lie algebra.

\section{Conjugacy classes in solvable Lie algebras}

Recall that given an element $X$, we say that
$X$ is of shape ${\mathcal I}$ if ${\mathcal I}$ is the ideal of smallest dimension containing $X$.

Factorizations of the adjoint group, related to normalizers of a given subalgebra, are very 
useful in computing the conjugacy class of the subalgebra under the adjoint group. It seems that 
for solvable algebras, the adjoint group can be factored in any order for a given basis of the 
algebra. However we have no proof of this statement. This explains the reason for the 
factorizations given below, which always exist and are adapted to a given element.

\textbf{Step 1.} Find all the ideals of $g$.

\textbf{Step 2.} Determine the possible
shapes of elements of $g$, starting from ideals of lowest dimension. For $X$ in $g'\,=\,
[g,\, g]$ embed $\langle
X\rangle $ a series of normalizers $$\langle X\,=\,X_{1}\rangle \,\subset\, \langle X_1,\, X_2\rangle
\,\subset\, \cdots \,\subset\, \langle X_1,\, X_2, \,\cdots,\, X_n
\rangle\, =\,g.$$ This gives a factorization of the adjoint group $e^{t_{n}{\rm ad}(X_{n})}
\cdots e^{t_{1}{\rm ad}(X_{1})}$.
If $X$ is not in $g'$, get a
factorization $e^{t_{1}{\rm ad}(X_{1})}\cdots e^{t_{r}{\rm ad}(X_{r})}e^{t_{r+1}{\rm ad}(X_{r+1})}
\cdots e^{t_{n}{\rm ad}(X_{n})}$ of the
adjoint group where $X_1 \,=\, X,\, \cdots ,\, X_{r}$ are independent mod $g'$ and $X_{r+1},\, \cdots,\, X_{n}$ are in $g'$.

\textbf{Step 3.} Use the factorizations in Step 2, to compute conjugacy classes of $1$-dimensional subalgebras of $g$: see the examples given below for this step.

\textbf{Step 4.} Any solvable subalgebra $h$ of $g$ can be obtained by a sequence of subalgebras
$$0\,=\,h_{0}\,\subset\, h_{1}\, \subset\, \cdots\, h_{d}\,=\,h ,$$ where
$\dim(h_{i+1}/h_{i})\,=\,1$ and $h_{i}$ is an ideal of $h_{i+1}$. Assume inductively that $\langle S_{1}\rangle ,\, \cdots,\,
\langle S_{k}\rangle$ are representatives of conjugacy classes of $d$-dimensional subalgebras. If $S$ is in this list and $N(S)\,=\,S$, then
$S$ cannot be extended to a higher dimensional algebra.

If $N(S)$ contains $S$ properly, then conjugacy classes of $1$-dimensional subalgebras of $N(S)/S$ determine all the $1$-dimensional
extensions of $S$.

Repeating this for all $S_{i}$, $1 \,\leq\, i \,\leq\, r$, gives all $(d+1)$-dimensional algebras, up to non-conjugacy. Compute the action of
the adjoint group --- adapted to a basis of $S_{i}$ --- to remove conjugates of $(d+1)$-dimensional algebras.

We conclude this section by stating a general result.

\begin{proposition}
 Let $g$ be a solvable Lie algebra with $g\,=\,
T+N$, where $N$ is a nilpotent ideal and $T$ is an abelian ad-digonalisable subalgebra. Then
every element $X\,=\,t+n$ with $t\,\in\, T$,\, $n\,\in\, N$ is conjugate to $\widetilde{X}
\,=\, t+ \widetilde{n}$ with $[ t,\, \widetilde{n}]\, =\,0$.
\end{proposition}

\section{Basic examples}

In this section we apply the algorithm to low dimensional examples, where computer assisted calculations can be avoided.

The following two lemmas have all the features of higher dimensional examples. They will be used repeatedly for computing conjugacy classes
in higher dimensions.

\begin{lemma}
Let $g$ be the $2$-dimensional algebra defined by the relation $[X,\,Y]\, =\, Y$. Then, the $1$-dimensional subspaces of $g$ are represented
by $\langle X\rangle $ and $\langle Y\rangle $.
\end{lemma}

\begin{proof}
As $Z(g)\,=\,0$ and $g'\,=\,\langle Y\rangle $, we see that $g$ has only one proper ideal, namely, $\langle Y\rangle $. Thus every $1$-dimensional
subspace is conjugate to $\langle Y\rangle $ or $\langle X+kY\rangle $. The adjoint group of $g$ factorizes as
$e^{{\rm ad} \langle Y\rangle } e^{{\rm ad}\langle X+kY\rangle }$. Therefore in computing conjugates of $\langle X+kY\rangle $ we
need only to consider its conjugate under $e^{t {\rm ad} Y}$.

Now ${\rm ad}Y(X+kY)\,=\,-Y$. Thus, $e^{t {\rm ad} Y}(X+kY)\,=\, X+kY-t Y$. Therefore taking $t\,=\,k$ we get
$\langle X+kY\rangle\, \sim\, \langle X\rangle $.
Hence representatives of $1$-dimensional subalgebras of $g$ are $\langle X\rangle $ and $\langle Y\rangle $.
\end{proof}

It is well known that a $3$-dimensional solvable Lie algebra $g$ is determined by the eigenvalues of $g/g'$ on $g'$.

\begin{lemma}
Let $g$ be a $3$-dimensional algebra with relations $[X,\,Y]\,=\,\lambda\cdot Y$,\, $[X,\,Z]\,=\, \mu\cdot Z$, where $\lambda \,\neq\, \mu$
and $\lambda \mu \,\neq\, 0$. The proper ideals of $g$ are $\langle Y\rangle $,\, $\langle Z\rangle $ and $\langle Y,\,Z\rangle $. The
conjugacy classes of $1$-dimensional subalgebras are $\langle X\rangle $, $\langle Y\rangle $, $\langle Z\rangle $, $\langle Y+Z\rangle $,
$\langle Y-Z\rangle $. The conjugacy classes of $2$-dimensional
subalgebras are represented by $\langle X,\,Y\rangle $,\, $\langle X, \,Z \rangle$,\, $\langle Y,\,Z\rangle$.
\end{lemma}

\begin{proof}
We have $Z(g)\,=\,0$ and $g'\,=\, \langle Y,\,Z\rangle$. The eigenspaces of $g/g'$ on $Z(g')=g'$ are $\langle Y\rangle $ and $\langle Z\rangle $.
Therefore, these are the only $1$-dimensional ideals of $g$.

The $2$-dimensional ideals of $g$ are given by the $1$-dimensional ideals of $g/\langle Y\rangle $ and $g/\langle Z\rangle .$ Now $g/\langle
Y\rangle \,=\, \langle \widetilde{X},\, \widetilde{Z}\rangle $ and $[\widetilde{X},\, \widetilde{Z}]\,=\,\lambda \widetilde{Z}$. Thus by Lemma
3.1, there is only one $1$-dimensional ideal of $g/\langle Y\rangle $, namely $\langle \widetilde{Z}\rangle $. Similarly, $g/\langle Z\rangle
$ has only one $1$-dimensional ideal namely $\langle \widetilde{Y}\rangle $.
Consequently, there is only one $2$-dimensional ideal of $g$, namely
$\langle Y,Z\rangle $.

Thus the proper ideals of $g$ are: $\langle Y\rangle $, $\langle Z\rangle $ and $\langle Y,Z\rangle $. Therefore, the
$1$-dimensional subspaces of $g$ are of the form $\langle Y\rangle$,\, $\langle Z\rangle$,\, $\langle Y+kZ\rangle ,\, \langle X+lY+mZ\rangle $.

In computing conjugates of $\langle Y+kZ\rangle $, we use the factorization of the adjoint group: $e^{\langle X\rangle }e^{\langle Y\rangle
}e^{\langle Y+kZ\rangle }$. As $Y$ and $Z$ commute in fact we need to consider conjugates only under $e^{{\rm ad}\langle X\rangle }$. Now
$[X,\,Y+kZ]\,=\,\lambda Y +k\mu Z$. Hence $e^{t ad(X)}(Y+kZ) \,=\, e^{t\lambda } Y + e^{t\mu} k Z$. Thus
$$\langle Y+kZ\rangle \,\sim\, \langle e^{t \lambda}Y
+ e^{t \mu k}k Z\rangle \,=\,\langle Y+e^{t(\mu-\lambda)}kZ\rangle .$$ Consequently, if $k\,\neq\, 0$, we have
$\langle Y+kZ\rangle \,\sim\, \langle Y+\epsilon Z\rangle$, where $\epsilon^2 \,=\,1$.

For computing conjugates of $\langle X+\lambda Y+\mu Z\rangle$, we use the factorization of the adjoint group: ${\rm ad}(g) \,=\, e^{\langle Z\rangle
}e^{\langle Y\rangle }e^{X+\lambda Y+\mu Z}$. In general, for a nilpotent Lie algebra, all of its elements operate as nilpotent
transformations on the Lie algebra. Now $$[Y,\, X+\lambda Y+\mu Z]\,=\, -\lambda Y,$$ Thus, $e^{{\rm ad}tY}(X+\lambda Y+\mu Z)
\,=\,X+\lambda Y+\mu Z-t\lambda
Y$. Therefore, $[X+\lambda Y+\mu Z] \,\sim\, [X+\mu Z]$, and conjugating $[X+\mu Z]$ by $e^{tZ}$ shows that $[X+\mu Z]\, \sim\, [X]$. Consequently,
representatives of conjugacy classes of $1$-dimensional subalgebras of $g$ are $\langle Y\rangle$,\, $\langle Z\rangle$,\,
$\langle Y+Z\rangle$,\, $\langle Y-Z\rangle $ and $\langle X\rangle$.

\textbf{$2$-dimensional subalgebras of $g$}:\,
A $2$-dimensional subalgebra of $g$ is obtained by extending a $1$-dimensional subalgebra. We have $N(\langle X\rangle )\,=\,\langle X\rangle $, so
$\langle X\rangle$ can't be extended. Now $N\langle Y\rangle \,=\,\langle Y,\,Z,\,X\rangle $. Thus $N\langle Y\rangle /\langle Y\rangle
\,=\, \langle
\widetilde{Z}, \widetilde{X}\rangle $. We have $[\widetilde{X},\, \widetilde{Z}]\,=\, \mu\cdot \widetilde{Z}$. By Lemma 3.1, the $1$-dimensional
subalgebras of $N\langle Y\rangle /\langle Y\rangle $ are $\langle \widetilde{X}\rangle $ and $\langle \widetilde{Z}\rangle $. Thus the
$2$-dimensional subalgebras obtained by extending $\langle Y\rangle $ are $\langle X,\,Y\rangle $ and $\langle Z,\,Y\rangle$. Similarly, the
two subalgebras obtained by extending $\langle Z\rangle $ are $\langle X,\,Z\rangle $ and $\langle Y,\,Z\rangle $. Now $N\langle Y+Z\rangle $ is
$\langle Y+Z,\,Z\rangle $ and $N\langle Y+Z\rangle /\langle Y+Z\rangle \,= \,\langle \widetilde{Z}\rangle $. Thus the $2$-dimensional subalgebra
with base $Y+Z$ is $\langle Y+Z,\,Z\rangle \,=\,\langle Y,\,Z\rangle $. Similarly, the $2$-dimensional algebra with base $Y-Z$ is $\langle Y-Z,
\,Z\rangle\,=\, \langle Y,\,Z\rangle $. Therefore, the representatives of conjugacy classes of $2$-dimensional subalgebras of $g$
are $\langle X,\,Y\rangle $,\, $\langle X,\,Z\rangle$ and $\langle Y,\,Z\rangle$.
\end{proof}

\subsection*{Example 4.1: subalgebras of a $3$-d solvable algebra that has no $1$-d ideal}
Let $g$ be a Lie algebra with basis $X\,=\,\partial_{x}$,\, $Y\,=\,\sin{x} \partial_{y}$,\, $Z\,=\,\cos{x}\partial_{y}$.

\textbf{Ideals:} In this case one can show directly that $g$ has no $1$-dimensional ideals and $g'$ is the only $2$-dimensional ideal of $g$.
Let us determine all the ideals and subalgebras of $g$ using the algorithms.

The algebra $g^{\mathbb{C}}$ is of the type considered in Lemma 3.1. It has a basis
$$\langle \partial_{x},\,\, e^{\sqrt{-1}x} \partial_{y},\,\, e^{-\sqrt{-1}x} \partial_{y}\rangle$$
as the eigenvalues of $\partial_{x}$ on $(g^{\mathbb{C}})'$ are $\sqrt{-1},\, -\sqrt{-1}$. Thus the ideals of $g^{\mathbb{C}}$ are
$\langle e^{\sqrt{-1}x} \partial_{y}\rangle $,\, $\langle e^{-\sqrt{-1}x} \partial_{y}\rangle $ and
$\langle e^{\sqrt{-1}x} \partial_{y}, \, e^{-\sqrt{-1}x} \partial_{y}\rangle $
by Lemma 3.1. Let $I$ be a proper ideal and $I^{\mathbb{C}}$ its complexification. If $\dim_{\mathbb{C}}I^{\mathbb{C}}\,=\,1$, we may
suppose that $I^{\mathbb{C}}\,=\,\langle e^{\sqrt{-1}x}\partial_{y}\rangle$. As $Re(a+\sqrt{-1}b)e^{\sqrt{-1}x}\partial_{y}
\,=\,a \cos{x} \partial_{y} + b \sin{x} \partial_{y}$, we see that $I^{\mathbb{C}}\cap g \,=\, I\, = \,
\langle \cos{x} \partial_{y},\, \sin{x} \partial_{y}\rangle $.

If $\dim_{\mathbb{C}}I^{\mathbb{C}}\, =\,2$, then $I^{\mathbb{C}}\,=\,\langle e^{\sqrt{-1}x}\partial_{y},\, e^{-\sqrt{-1}x}\partial_{y}\rangle $.
Now $I^{\mathbb{C}}\cap g $ is given by the real part of $z(e^{\sqrt{-1}x}\partial_{y}) + w (-e^{\sqrt{-1}x}\partial_y)$, where
$z,\,w \,\in\, \mathbb{C}$. The real part works out to
be $\langle \cos{x}\partial_{y},\, \sin{x} \partial_{y}\rangle$. Thus the only proper ideal of $g$ is $g'$.

\textbf{Subalgebras:} There is only one proper ideal of $g$, namely $\langle Y,\,Z\rangle$. Thus the $1$-dimensional subalgebra is a
subalgebra of $\langle Y,\,Z\rangle $ or it is of the form $\langle X+ a Y +bZ\rangle$. To consider conjugacy classes of $1$-dimensional
subalgebras of the ideal $\langle Y,\,Z\rangle $, we only need to consider its conjugate under the $1$ parameter subgroup generated by
$\partial_{x}$ because $[Y,\,Z]\,=\,0$. Now $e^{{\rm ad} tX}$ operates as rotations on the plane $\langle Y,\,Z\rangle$. Thus every $1$-dimensional
subspace of $\langle Y,\,Z\rangle$ is conjugate to $\langle Y\rangle$.

If we have a $1$-dimensional subalgebra $\langle \widetilde{X} \,=\, X+aY+bZ\rangle$, then the adjoint group factorizes as $e^{\langle
{\rm ad}(Z)\rangle } e^{\langle {\rm ad}(Y)\rangle } e^{\langle{\rm ad}(\widetilde{X})\rangle }$. Thus to compute conjugates of $\langle \widetilde{X}\rangle $ we
need only compute its conjugates under $e^{{\rm ad}\langle Z\rangle } e^{{\rm ad}\langle Y\rangle }$. Both ${\rm ad}(Z)$ and
${\rm ad}(Y)$ must be nilpotent. Now ${\rm ad}(Y)(X+aY+bZ)\,=\,Z$ and ${\rm ad}(Y)(Z)\,=\,0$. Thus, $e^{t{\rm ad}(Y)}(X+aY+bZ)
\,=\, X+aY+bZ+tZ$. Setting $t\,=\,-b$, we see that $\langle X+aY+bZ\rangle \,\sim\, \langle
X+aY\rangle $. Moreover, ${\rm ad}(Z)(X+aY)\,=\,-Y$ and ${\rm ad}(Z)(Y)\,=\,0$. Hence, $e^{s\cdot {\rm ad} Z}(X+aY)\,=\,X+aY-sY$.
Setting $s\,=\,a$, we see that $\langle X+aY\rangle \,\sim\,
\langle X\rangle $. Hence, representatives of $1$-dimensional subalgebras of $g$ are $\langle X\rangle $ and $\langle Y\rangle $. Any
$2$-dimensional subalgebra is an extension of a $1$-dimensional subalgebra. As $N\langle X\rangle \,=\, \langle X\rangle $,
it follows that $\langle X\rangle $
cannot be extended. As $N\langle Y\rangle \,=\,\langle Y,\,Z\rangle $, we have $N\langle Y\rangle /\langle Y\rangle \,=\,\langle Z\rangle$,
and hence the $1$-dimensional extension of $Y$ is $\langle Y,\,Z\rangle \,=\,g'$.

\begin{remark}
Let ${\mathfrak g}\,=\, \langle X,\, Y,\, Z\rangle$ with
$$
[X,\, Y]\,=\, aY+bZ,\,\ \ [X,\, Z]\,=\, -bY+aZ
$$
such that $a^2+b^2\, \not=\, 0$. One can check that $Z(g)\,=\,0$ and ${\rm ad}(X)$
restricted to $Z(g')$ has no real eigenvalues. Thus $g$ has no $1$--d ideal and $1$--d algebra
are conjugate to $X$,\, $Y$ and any $2$--d subalgebra coincides with $g'\,=\,\langle Y,\,Z\rangle$.
\end{remark}

\subsection*{Example 4.2: subalgebras of the unique nilpotent nonabelian $3$-d algebra}

Let $g$ be the Lie algebra with basis $X,\,Y,\,Z$ and non-zero commutator $[X,\,Y]\,=\,Z$.

\textbf{Ideals:} Now $g'\,=\,\langle Z\rangle \,=\,Z(g)$. Thus $g$ has only one $1$-dimensional ideal, namely $\langle Z\rangle $. Now $g/\langle
Z(g)\rangle $ is abelian, so its ideals are $\langle a \widetilde{X}+ b\widetilde{Y}\rangle $, and
hence are $\langle \widetilde{Y}\rangle$, $\langle
\widetilde{X}+k\widetilde{Y}\rangle $. Therefore, the $2$-dimensional ideals are $\langle Y,\,Z\rangle $ and $\langle X+kY,\,Z\rangle $. Thus all
the ideals of $g$ are $\langle Z\rangle $,\, $\langle Y,\,Z\rangle $ and $\langle X+kY,\,Z\rangle $.

\textbf{Subalgebras:} From the description of the ideals, every $1$-dimensional subspace is of the form $\langle Z\rangle $,\, $\langle
Y+mZ\rangle $ or $\langle X+k Y + lZ\rangle$. Consider a $1$-dimensional subspace of the form $\langle X+kY+lZ\rangle \,=\,\langle
\widetilde{X}\rangle$. A basis of $g$ is $\widetilde{X},\, Y,\,Z$. As $\widetilde{X}$,\, $Z$ centralize $\widetilde{X}$, we need to
consider only the $1$-parameter subgroup generated by $Y$ to determine conjugates of $\langle \widetilde{X}\rangle $.
We have ${\rm ad}(Y)(\widetilde{X})\,=\,Z$ and ${\rm ad}(Y)(Z)\,=\,0$. Thus,
$e^{t{\rm ad}(Y)}(\widetilde{X}) \,=\, \widetilde{X}+tZ$. Hence, $\langle X+kY+lZ\rangle \,\sim\, \langle X+k Y\rangle$. Similarly,
$\langle Y+mZ\rangle \,\sim\, \langle Y\rangle$. Therefore, representatives of conjugacy classes of $1$-dimensional subalgebras are
$\langle Z\rangle $,\, $\langle
Y\rangle $ and $\langle X+kY\rangle$. The $2$-dimensional subalgebras are obtained by extending these representatives. Now $N(Z)\,=\,\langle
X,\,Y,\,Z\rangle $. Thus, $N(Z)/\langle Z\rangle \,=\,\langle \widetilde{X},\, \widetilde{Y}\rangle$ which is abelian. Hence, the $1$-dimensional
subalgebras are $\langle \widetilde{Y}\rangle $ and $\langle \widetilde{X}+k \widetilde{Y}\rangle $. Thus extensions of $\langle Z\rangle $ are
$\langle Z,\,Y\rangle $ and $\langle Z,\,X+kY\rangle $. Both are ideals of $g$. The $2$-dimensional algebras obtained by extending
$\langle Y\rangle$ come from classes of $1$-dimensional subalgebras of $N\langle Y\rangle /\langle Y\rangle
\,=\, \langle Y,Z\rangle /\langle Y\rangle $. Hence, $\langle Y,\,Z\rangle $ is the only extension of $\langle Y\rangle $. Similarly,
$N\langle X+kY\rangle \,=\, \langle X+kY,\,Z\rangle $. Thus the only
$2$-dimensional extension of $\langle X+kY\rangle $ is $\langle X+kY,\,Z\rangle$. Hence the $2$-dimensional subalgebras of $g$ are $\langle
Y,\,Z\rangle$,\, $\langle Z,\,X+kY\rangle $, where $k\,\in \,\mathbb{R}$. Therefore all $2$-dimensional subalgebras of $g$ are ideals of $g$.

\subsection*{Example 4.3: subalgebras of a $4$-d algebra from physics}

This is a well-known example considered in the literature (see \cite[p.~203]{Ol2}).
Let $g\,=\,\langle e_{1},\, e_{2},\, e_{3}, \,e_{4}\rangle$, where $e_{1} \,=\, t\partial_{x}+\partial_{u}$,\, $e_{2}
\,=\, x\partial_{x}+3t\partial_{t}-2u\partial_{u}$,\, $e_{3}\,=\, \partial_{t}$ and $e_{4}\,=\,\partial_{x}$.

\textbf{Ideals:} The commutators of the basis are: $[e_{1},\,e_{2}]\,=\,-2e_{1}$,\, $[e_{1},\,e_{3}]\,=\,-e_{4}$,\,
$[e_{2},\,e_{3}]\,=\,-3e_{3}$ and $[e_{2},\,e_{4}]\,=\,-e_{4}$. The computations here are very similar to those in Lemma 3.1 and
Lemma 3.2. We omit
most details. We have $Z(g)\,=\,0$ and $Z(g')\,=\,\langle e_{4}\rangle $. Working in $g^{\mathbb{C}}$, we find that there is only one
$1$-dimensional ideal of $g^{\mathbb{C}}$, namely $\langle e_{4}\rangle $. Working in quotients of $g^{\mathbb{C}}$ by ideals already
obtained, we see that all the proper ideals of $g^{\mathbb{C}}$ are $\langle e_{4}\rangle $,\, $\langle e_1,\, e_4\rangle $,\, $\langle e_3,\,
e_4\rangle $ and $\langle e_1,\, e_3, \,e_4\rangle$. As these ideals have generators in $g$, we see that the above list is also a list of all
proper ideals of $g$. Therefore, the $1$-dimensional subalgebras of $g$ are $\langle e_4\rangle $,\, $\langle e_1+a e_4\rangle $,\, $\langle e_3
+b e_4\rangle$,\, $\langle e_1+ce_3 +d e_4\rangle $ and $\langle e_2 +f e_1 + g e_3 + he_4\rangle$.

\textbf{Subalgebras:} To find the subalgebras we follow the same procedure as before which is based on Lemma 3.1 and Lemma 3.2. For any element $X$
of $g'$, we use a series of normalizers $$\langle X\,=\,X_1\rangle \,\subset\, \langle X_1,\, X_2\rangle \,\subset\,
\langle X_1, \,X_2,\, X_3\rangle \,\subset\, \langle X_1,\, X_2,\, X_3,\, X_4\rangle$$ to factorize the adjoint group. For
$X\,=\, e_2 + fe_1 + g e_3 +h e_4 $ we have $N(X)\,=\,X$, and $\langle X,\, e_1,\, e_3,\,
e_4\rangle $ is a basis of $g$ with $X$ normalizing the ideal $\langle e_1,\, e_3,\, e_4\rangle$.

In this case, the adjoint group factorizes as $e^{{\rm ad}\langle e_4\rangle } e^{{\rm ad}\langle e_3\rangle }
e^{{\rm ad}\langle e_1\rangle } e^{{\rm ad}\langle X\rangle } $. We use such factorizations to find conjugates of $\langle X\rangle $
under the adjoint group. We omit the details as this is similar to the computations done before. For the convenience of the
reader, we discuss conjugacy of two types of elements. Let $X\,=\,\langle e_1 + ce_3 + de_4\rangle $. Here we use the series of normalizers
$$\langle X\rangle \,\subset\, \langle X,\, e_4\rangle\,\subset \, \langle X,\,e_4,\, e_3\rangle
\,\subset\, \langle X,\, e_4,\, e_3,\, e_2\rangle .$$ The adjoint group thus factorizes as
$e^{{\rm ad}\langle e_2\rangle } e^{{\rm ad}\langle e_3\rangle } e^{{\rm ad}\langle e_4\rangle } e^{{\rm ad}\langle X\rangle }$. As
$X$ and $e_{4}$ commute, we only need to compute conjugates under
$e^{s {\rm ad}\langle e_4\rangle } e^{t {\rm ad}\langle e_3\rangle }$. Both of ${\rm ad} e_{3}$ and ${\rm ad} e_{2}$ must be nilpotent
as both are in $g'$. We find that $$\langle X\rangle \, \sim\, \langle e_{1}+c e_{3}\rangle \,\sim\, \langle e^{2s}e_{1}
+c e^{-3s}e_{2}\rangle\,\sim\, \langle e_{1} +c e^{-5s}e_{2} \rangle\,=\, \langle e_{1}+ \epsilon e_2\rangle ,$$ where $\epsilon
\,=\, \pm 1$. For $X\,=\, e_2 + f e_1 + g e_3 +h e_4$ we use the factorization $e^{{\rm ad}\langle e_4\rangle } e^{{\rm ad}\langle e_3\rangle }
e^{{\rm ad}\langle e_1\rangle } e^{{\rm ad}\langle X\rangle }$. We can ignore the last factor because it stabilizes $\langle X\rangle$.

By computations similar to the types already considered, we find that $\langle X\rangle\,\sim\, \langle e_{2}\rangle $. Thus
representatives for conjugacy classes of $1$-dimensional subalgebras are $\langle e_4\rangle $,\, $\langle e_1\rangle $,\, $\langle e_3\rangle $,
\,$\langle e_1+e_3\rangle $,\, $\langle e_1-e_3\rangle $ and $\langle e_2\rangle$.

\textbf{$2$-dimensional Subalgebras of $g$:}\, Two dimensional algebras are obtained by extending $1$-dimensional algebras. Thus we need to find
normalizers of representatives of $1$-dimensional algebras, say $H_{i}$, and find representatives of $1$-dimensional conjugacy classes in
$N(H_{i})/H_{i}$. There will be algebras of dimension at most 3 and their subalgebras have already been determined in Lemma 3.1 and Lemma 3.2.
We will just list the classes of $2$-dimensional algebras. However, we indicate briefly how all the $2$-dimensional subalgebras which are
extensions of $\langle e_{4}\rangle$ can be obtained.

Now $N(e_{4})\,=\,g$ and therefore $N(e_{4})/\langle e_{4}\rangle \,=\,\langle \widetilde{e}_{2}, \,\widetilde{e}_{1},
\widetilde{e}_{3}\rangle$. As $[\widetilde{e}_{2},\,\widetilde{e}_{1}] \,=\, 2 \widetilde{e}_1$,\, $[\widetilde{e}_{2},\, \widetilde{e}_{3}]\,
= \,-3 \widetilde{e}_{3}$, Lemma 3.2 gives the conjugacy classes of $1$-dimensional subalgebras of
$N\langle e_{4}\rangle /\langle e_{4}\rangle$ as
$\langle \widetilde{e}_2\rangle $,\, $\langle\widetilde{e}_1\rangle $,\, $\langle \widetilde{e}_3\rangle $,\, $\langle \widetilde{e}_1 +
\widetilde{e}_{3}\rangle $ and $\langle \widetilde{e}_1 - \widetilde{e}_{3}\rangle $. Hence, the $2$-dimensional subalgebras obtained by
extending $\langle e_{4}\rangle$ are $\langle e_4,\, e_2\rangle $,\, $\langle e_4,\, e_1\rangle $,\, $\langle e_4,\, e_3\rangle $,\,
$\langle e_4, \,e_1+e_3\rangle $ and $\langle e_4, \,e_1-e_3\rangle$.

The other cases are similar and easier. One uses Lemma 3.1 for these cases. After removing repetitions, we obtain the following list of
$2$-dimensional subalgebras: $\langle e_4,\, e_2\rangle $,\, $\langle e_4,\, e_1\rangle $,\, $\langle e_4,\, e_3\rangle $,\,
$\langle e_4,\, e_1+e_3\rangle$,\, $\langle e_4,\, e_1-e_3\rangle $ and $\langle e_2, \,e_3\rangle $. In computing conjugates of these subalgebras, we can choose bases
complementary to the basis of a given subalgebra to simplify the computations. Also, we need to do this only for subalgebras that are not
ideals.

The conjugates of $[e_{4} \wedge e_2]$ are $$e^{s {\rm ad}(e_{1})}\, e^{t {\rm ad} (e_{3})} [e_4 \wedge e_2]\,=\, [e_4 \wedge e_2 + 3t e_4
\wedge e_3 -2s e_4 \wedge e_1].$$ Thus $[e_4 \wedge e_2]$ is not conjugate to any other $2$-dimensional subalgebra listed above. By similar
computations, we find that all the $2$-dimensional algebras listed above are non-conjugate.

Hence representative for classes of $2$-dimensional subalgebras of $g$ are $\langle e_4,\, e_2 \rangle $,\, $\langle e_4,\, e_1\rangle $,\,
$\langle e_4,\, e_3\rangle $,\, $\langle e_4,\, e_1+e_3\rangle $,\, $\langle e_4,\, e_1-e_3\rangle $ and $\langle e_2,\, e_3\rangle$.

\textbf{$3$-dimensional Subalgebras of $g$:}\, To obtain $3$-dimensional algebras we need to compute $1$-dimensional extensions of
representatives of conjugacy classes of $2$-dimensional algebras. If $A$ is such an algebra, $N(A)/A$ is of dimension at most two and Lemma 3.1
gives the $1$-dimensional extensions in case $N(A)/A$ is non-abelian. Doing this for all classes of $2$-dimensional algebras we get the
following extensions: $\langle e_4,\, e_1,\, e_2\rangle$,\, $\langle e_4,\, e_1,\, e_3\rangle$ and
$\langle e_4,\, e_3,\, e_2\rangle$. These algebras are non-conjugate because
\begin{align}
e^{t {\rm ad} (e_{3})}(e_{4} \wedge e_1 \wedge e_2) &\,=\, e_4 \wedge (e_1 -t e_4) \wedge (e_2 + 3te_3)\nonumber \\
&=\,e_4\wedge e_1 \wedge e_2 + 3t e_4 \wedge e_1 \wedge e_2.
\end{align}
Thus $[e_4 \wedge e_1 \wedge e_2]$ is not conjugate to $[e_4\wedge e_3 \wedge e_2]$. Moreover, $e^{t {\rm ad} (e_{2})}
[e_{4} \wedge e_1 \wedge e_3] \,=\, [e_4\wedge e_1 \wedge e_3]$.

\subsection*{Example 4.4: subalgebras of a $6$-d algebra from physics}

In this example, we indicate briefly how the subalgebras of non-solvable algebra can be determined up to conjugacy.

We take for $g$ the $6$ dimensional algebra given in \cite{AO}. Writing the basis elements $\Xi_{i}$ as $e_{i}$, $1\,\leq\, i \,\leq\, 6$, and using
standard Maple commands, we find that the Levi decomposition is
\begin{equation}
\langle e_{1},\, e_{3},\, e_{5} \rangle \oplus \langle e_{2},\, e_{3}- 2e_{4},\, e_{6} \rangle\, =\, r \oplus s,
\end{equation}
where $r\, \subset\, g$ is the nilpotent radical, and $s$ is the semisimple part.
We have $[e_1,\, e_5]\,=\,-e_{3}$ so the radical is the algebra discussed in Example 4.2.

Now $e_{2},\, e_{6}$ are eigenvectors of ${\rm ad}(e_3-2e_{4})$ with opposite eigenvalues and $[e_2,\, e_{6}]\,=\,-2(e_3 -2e_{4})$. Thus
by scaling the
generators of semisimple part, we see that it has generators $X,\,Y,\,H$ with $[H,\,X]\,=\,2X$,\, $[H,\,Y]\,=\,-2Y$ and $[X,\,Y]
\,=\,H$. Hence the semisimple part is
$\mathfrak{sl}(2, \mathbb{R})$. Its maximal solvable subalgebras are $\langle X-Y \rangle$ and $\langle H,\,X \rangle$. Therefore, the maximal
solvable subalgebras of $g$ are $\langle X-Y \rangle +r $ and $\langle H,X \rangle +r$. The first algebra, except for a central element is the
algebra in Example 4.1. The second algebra is similar to Example 4.2. Thus we know, up to conjugacy in the maximal solvable subalgebras, all
the solvable algebras. The adjoint group of $g$ is $\langle e^{{\rm ad}(\mathfrak{sl}(2,\mathbb{R})) }\rangle\,e^{{\rm ad}(r)}$.

Using the decomposition $\langle e^{{\rm ad}(\mathfrak{sl}(2,\mathbb{R})) }\rangle\,=\,
e^{\langle {\rm ad}(X-Y)\rangle} e^{\langle {\rm ad}(H)\rangle}e^{\langle {\rm ad}(X)\rangle}$, we can determine all classes of solvable
subalgebras.

The non-solvable algebras with a proper Levi decomposition must arise as $N(U)$, where $U$ is 
solvable. Finally $\mathfrak{sl}(2,\mathbb{R})$ is the only semisimple algebra of $g$ --- up to 
conjugacy in the adjoint group. We leave further details to the reader.

\section*{Acknowledgements}

We thank the referee for helpful comments.

\end{document}